\documentclass[11pt,3p]{elsarticle}

\usepackage{natbib}
\usepackage{amsmath,amssymb,mathrsfs}
\usepackage{color}
\usepackage{graphicx}
\usepackage{subfigure}
\usepackage{soul}
\usepackage{bbding}
\usepackage{lineno,hyperref}
\usepackage{multirow, booktabs}
\modulolinenumbers[5]

\def\bu{\mathbf{u}}

\def\eps{\varepsilon}

\def\bx{\mathbf{x}}

\def\grad{\nabla}
\def\laplace{\Delta}

\begin{document}

\begin{frontmatter}

\title{
A Unified Weighted-Loss Physics-Informed Neural Network for Boundary Layer Problems in Singularly Perturbed PDEs}

\author[NCU,NCTS]{Wei-Fan Hu} 
\ead{wfhu@math.ncu.edu.tw}
\author[NCU]{Shi-Xiang Zhong}
\author[NCHU]{Po-Wen Hsieh}
\ead{pwhsieh@nchu.edu.tw}
\author[NYCU]{Chung-Kai Chen}
\author[NYCU,NCTS]{Te-Sheng Lin\corref{cor1}}
\ead{teshenglin@nycu.edu.tw}

\cortext[cor1]{Corresponding author}

\address[NCU]{Department of Mathematics, National Central University, Taoyuan 320317, Taiwan}
\address[NCHU]{Department of Applied Mathematics, National Chung Hsing University, Taichung 402202, Taiwan}
\address[NYCU]{Department of Applied Mathematics, National Yang Ming Chiao Tung University, Hsinchu 300093, Taiwan}
\address[NCTS]{National Center for Theoretical Sciences, National Taiwan University, Taipei 106216, Taiwan}

\begin{abstract}

Singularly perturbed partial differential equations arise in many applications, including magnetohydrodynamic duct flows, chemical reaction–transport systems, and Poisson–Boltzmann electrostatics. These problems are characterized by sharp boundary layers and pronounced multiscale behavior, posing significant challenges for numerical methods. Existing approaches, particularly machine learning based methods, often rely on explicit asymptotic decompositions or specialized architectures, increasing implementation complexity and leading to optimization imbalance in stiff regimes. In this work, we propose a unified learning framework based on a weighted loss formulation within the standard physics--informed neural network setting. The proposed method requires only prior knowledge of the boundary layer thickness, while the boundary layer locations are automatically identified during training. The resulting formulation avoids problem specific architectural modifications and remains applicable across different equation types. Numerical experiments on both scalar and coupled reaction-diffusion and convection-diffusion-reaction systems, defined on regular and irregular domains, demonstrate robust performance for boundary layer thickness as small as $10^{-10}$ while maintaining high solution accuracy.
\end{abstract}
\end{frontmatter}


\section{Introduction}

Singularly perturbed problems arise in a broad range of scientific and engineering applications and serve as fundamental models for multiscale phenomena. Representative examples include turbulent wave-current interactions~\cite{ZRZZ22}, reaction-diffusion processes~\cite{M97}, optimal control problems~\cite{K84}, and magnetohydrodynamic flows~\cite{HY09}. The numerical approximation of such problems is well known to be challenging~\cite{RM15, MQQX25, ZSY25, AG22}. The primary difficulty lies in the presence of localized solution structures, such as boundary or interior layers, when the perturbation parameter tends to zero. These sharp gradients lead to significant loss of accuracy and numerical instability for standard discretization methods, even in one-dimensional settings~\cite{RM15, RST08}. In particular, when the mesh size is not sufficiently small relative to the perturbation parameter, classical schemes such as central difference methods or Galerkin finite element methods with piecewise linear basis functions can produce spurious oscillations. Such oscillations may appear not only within the layer regions but also contaminate the smooth part of the solution.

To address these difficulties, a variety of numerical strategies have been developed, including finite difference schemes with enhanced accuracy or stabilization~\cite{R94, TD07, HYY18}, layer-adapted mesh techniques~\cite{MOS02, Lin09, OS09}, stabilized or higher-order finite element methods~\cite{BH82, FC95, KVKS19}, and multiscale finite element approaches~\cite{BR94, BFR98, FRV05, FMTV05, HY10}. Despite these advances, many classical approaches rely heavily on \emph{a priori} information about the layer location and thickness to design fitted meshes or select stabilization parameters. Moreover, mesh-dependent methods are generally difficult to extend to higher-dimensional problems or complex geometries due to increased computational cost and implementation complexity. For coupled singularly perturbed systems, the interaction and overlap of multiple layers introduce additional numerical challenges, often resulting in reduced robustness and limited scalability~\cite{MS03}.

In recent years, neural-network-based methods have attracted increasing attention within the scientific community. Among these approaches, physics-informed neural networks (PINNs)~\cite{raissi2019} have emerged as a powerful framework for solving both forward and inverse problems, demonstrating promising performance across a wide range of scientific and engineering applications. The central idea of PINNs is to minimize a loss function associated with a given problem, which incorporates the residuals of the governing partial differential equation (PDE), together with the corresponding boundary and initial conditions. In contrast, deep operator networks (DeepONets)~\cite{LJPZK21} learn linear or nonlinear differential operators by approximating mappings between input data and solution fields, rather than enforcing the governing equations through a loss function. The success of both PINNs and DeepONets can be attributed to the expressive power of neural network representations~\cite{Cybenko1989, HSW89, CC95}, which enables them to approximate smooth solutions in appropriate function spaces.

Despite the success of deep learning-based approaches for PDEs with smooth solutions, both standard PINN and DeepONet frameworks often encounter convergence difficulties when solutions exhibit localized steep variations. Such behavior arises, for example, in boundary layers of singularly perturbed problems and, more broadly, in stiff and multiscale PDEs. These challenges underscore the need for more robust learning strategies tailored to stiff and multiscale problems. 

To address these issues, several approaches incorporate ideas from classical perturbation theory. The boundary-layer PINN (BL-PINN)~\cite{ACD23} mimics the construction of inner and outer asymptotic expansions and enforces a matching condition to obtain a uniformly valid solution. Numerical results demonstrate its effectiveness in capturing boundary-layer behavior. However, the method requires the construction of multiple neural networks to represent different components of the solution and introduces additional tuning parameters in the matching condition. Similarly, the scale-decomposed PINN (SD-PINN)~\cite{ZH26} decomposes the solution into inner, outer, and matching components without explicitly constructing asymptotic expansions. Nevertheless, its loss function involves multiple coupled residual terms, which may lead to slow training and optimization imbalance. Moreover, for some benchmark problems, the resulting accuracy shows limited improvement over leading-order asymptotic approximations.

Exploiting asymptotic decay structures, the approaches in~\cite{CGHJ24, WZH24} incorporate problem-specific decay functions into the neural network representation, improving accuracy in boundary-layer regions. However, these methods rely heavily on prior knowledge of layer locations and asymptotic structure. Beyond asymptotic-based strategies, hybrid approaches such as~\cite{YG24} combine neural networks with classical discretizations, for example by learning stabilization parameters in finite element methods. While effective for certain convection–diffusion problems, such methods increase model complexity, particularly when coefficient fields must be represented by convolutional networks, and their overall performance remains tied to the underlying numerical scheme and optimization accuracy.

To overcome these challenges, we propose a unified learning framework for a broad class of singularly perturbed PDEs, including scalar equations and coupled systems on both regular and irregular domains, covering linear and nonlinear formulations. In contrast to the asymptotic-based approaches, the proposed framework requires only prior knowledge of the boundary-layer thickness, rather than a full asymptotic characterization. Remarkably, the proposed network can autonomously detect boundary-layer locations during training, eliminating the need for \emph{a priori} information on layer locations. Although the solution representation also adopts a decomposition strategy, the key distinction lies in the weighted loss formulation, which is applicable to a broad range of PDEs with diverse underlying physical characteristics. Numerical experiments demonstrate that the framework remains stable and accurate even for extremely small perturbation parameters (down to $\eps = 10^{-10}$), achieving solution errors on the order of $10^{-7}$--$10^{-8}$. To the best of our knowledge, such performance for perturbation parameters of this magnitude has not been reported in existing neural-network-based solvers for singularly perturbed problems.

The remainder of this paper is organized as follows. Section~\ref{sec:formulation} introduces the problem formulation and the theoretical background. The proposed neural network decomposition strategy is presented in Section~\ref{sec:NN}, followed by the implementation details and the training framework in Section~\ref{sec:implementation}. Numerical experiments are reported in Section~\ref{sec:results}. Finally, concluding remarks are given in Section~\ref{sec:conclusion}.

\section{Problem formulation}\label{sec:formulation}

Given a small parameter $0 < \eps \ll 1$, we aim to develop a neural-network-based framework for solving a general class of $d$-dimensional ($d=1,2$) singularly perturbed elliptic problems defined on a domain $\Omega \subset \mathbb{R}^d$, subject to Dirichlet boundary conditions on $\partial \Omega$:
\begin{align}\label{Eq:PDE}
\begin{split}
\mathcal{L}_\eps(u) &= f, \quad \text{in } \Omega, \\
u &= g, \quad \text{on } \partial \Omega.
\end{split}
\end{align}
Here, $\mathcal{L}_\eps$ denotes a second-order elliptic operator with linear or nonlinear lower-order terms, in which the small parameter $\eps$ multiplies the highest-order derivative. This singular perturbation gives rise to rapid variations in the solution profile, typically manifested as boundary or interior layers.
We note that other types of boundary conditions, such as Neumann or Robin conditions, do not alter the main ingredients of the proposed framework.
Although~Eq.~(\ref{Eq:PDE}) is written in scalar form for conciseness, the proposed approach extends straightforwardly to systems of singularly perturbed problems, which will be addressed in subsequent sections.

It is well known that the solution $u$ of the boundary value problem~\eqref{Eq:PDE} may exhibit localized steep gradients in boundary or interior layers. Throughout this paper, we restrict our attention to problems involving boundary layers only and assume that the boundary layer width scales as $O(\eps)$.

Based on asymptotic expansion theory, the solution can be approximated by combining an outer approximation and an inner approximation. The outer approximation, obtained in the reduced problem ($\eps = 0$), describes the solution behavior away from the boundary layers, whereas the inner approximation, constructed through an appropriate spatial rescaling, captures the rapid variations within the narrow boundary layers. By imposing suitable matching conditions between these approximations, one obtains a uniformly valid approximation. Inspired by this classical theoretical framework, we assume that the solution $u$ of Eq.~(\ref{Eq:PDE}) admits the decomposition
\begin{align}\label{Eq:decomp}
u = u_r + u_s,
\end{align}
where $u_r$ denotes the \emph{regular} component characterizing the smooth behavior of the solution, and $u_s$ denotes the \emph{singular} component capturing the sharp variations within the boundary layers. Hereafter, the subscripts $r$ and $s$ are used to indicate the regular and singular components, respectively. These two components are constructed to reflect the behavior of the classical outer and inner approximations. It is noteworthy that both $u_r$ and $u_s$ are defined over the global domain $\Omega$, rather than being restricted to the regions in which the classical outer and inner approximations are conventionally valid. We remark that the decomposition in Eq.~(\ref{Eq:decomp}) is not unique, as the regular and singular components may differ by an additive constant while still yielding the same solution $u$.

\subsection{Model problems}\label{subsec:model}

To demonstrate the accuracy, generality, and robustness of the proposed approach, we consider several classes of singularly perturbed partial differential equations. All equations are defined on bounded domains with appropriate boundary conditions. Throughout this work, we focus on parameter regimes that induce boundary layer structures, while problems involving interior layers are not considered in the present study. The selected model problems encompass the following categories:
\begin{itemize}
\item \textbf{Reaction--diffusion equations.}
We consider equations of the form
\begin{align*}
-\eps^2 \laplace u + \gamma u = f,
\end{align*}
where $\gamma > 0$ is a given reaction coefficient and $f$ is a prescribed source term.

\item \textbf{Convection--diffusion--reaction equations.}
We investigate the following model: 
\begin{align*}
-\eps \laplace u + \mathbf{a}\cdot\nabla u + \gamma u = f,
\end{align*}
where $\mathbf{a}$ is a given velocity field with $\|\mathbf{a}\| = O(1)$ and $\gamma > 0$ denotes a reaction coefficient. The case $\gamma = 0$ reduces to the classical convection--diffusion equation.

\item \textbf{Poisson--Boltzmann equations.}
We consider the nondimensionalized nonlinear Poisson--Boltzmann equation
\begin{align*}
-\eps^2 \laplace u + \exp(u) - \exp(-u) = 0,
\end{align*}
which arises in electrostatic modeling of electrolytes. This equation is characterized by a strong coupling between the singular perturbation and the nonlinear reaction term, which significantly complicates the solution profile within boundary layers.

\item \textbf{Systems of convection--diffusion--reaction equations.}
We consider systems of the form
\begin{align*}
-\eps \laplace \bu + A \bu_x + B \bu_y + C \bu = \mathbf{f},
\end{align*}
where $\bu = [u_1, u_2, \ldots, u_n]^\top$ is a vector-valued unknown; $A$, $B$, and $C \in \mathbb{R}^{n \times n}$ are given coefficient matrices; and $\mathbf{f} = [f_1, f_2, \ldots, f_n]^\top$ denotes the source term. In the one-dimensional setting, the system reduces to the corresponding form involving only spatial derivatives in $x$.
\end{itemize}

Having established the mathematical formulations of these model problems, the subsequent section details the construction of the neural network architectures and the corresponding efficient training strategies.

\section{Methodology: solution decomposition via neural networks}\label{sec:NN}
For small values of $\eps$, the standard PINN formulation frequently fails to achieve convergence when solving Eq.~(\ref{Eq:PDE}) due to the singular nature of the solution, characterized by stiff gradients and rapidly varying higher-order derivatives within the boundary-layer regions. These features pose significant challenges for the optimization landscape and the approximation capacity of standard networks. We elucidate the mathematical formulation of the solution decomposition in Eq.~(\ref{Eq:decomp}) and detail its consistent extension from one-dimensional to two-dimensional domains.

\subsection{One-dimensional case}\label{subsec:1D}
We begin by introducing the fundamental concept of the solution decomposition in Eq.~(\ref{Eq:decomp}) in the one-dimensional setting, with $\Omega = [a,b]$. 
It is emphasized that the regular component $u_r$ captures the smooth variations of the solution, conceptually aligning with the classical outer approximation. Given its high regularity, this component can be efficiently approximated by a conventional multilayer perceptron (MLP).
By construction, $u_r$ and its derivatives up to second order remain of magnitude $O(1)$ over the entire domain.

On the other hand, the singular component $u_s$ and its derivatives remain nearly constant outside the boundary-layer regions, while scaling arguments indicate that its first derivative is of order $O(1/\eps)$ since the boundary-layer width is $O(\eps)$. This behavior can be captured by combining a continuous sigmoidal-type function with a level-set function. Specifically, for a layer near $x_0$, we consider the  sigmoidal activation function $\sigma(x)=1/(1 + e^{-x})$. Then
\begin{equation}
\sigma\Big(\frac{x-x_0}{\eps}\Big)=\frac{1}{1+e^{-(x-x_0)/\eps}}.
\end{equation}
This function exhibits a rapid transition near $x_0$. In particular,
\begin{equation}
\left(\sigma\Big(\frac{x-x_0}{\eps}\Big)\right)'
= \sigma\Big(\frac{x-x_0}{\eps}\Big)
  \left(1-\sigma\Big(\frac{x-x_0}{\eps}\Big)\right)\frac{1}{\eps},
\end{equation}
indicates that the derivative is of order $O(1/\eps)$ within a region of width $O(\eps)$, while remaining nearly zero elsewhere due to the exponential decay of the sigmoidal function away from the layer. These properties are consistent with the characteristic behavior of the singular component. Consequently, motivated by this observation, we assume that the singular component can be represented as
\begin{align}\label{Eq:1D_singular}
u_s(x) = u_{s,L}\Big(\frac{\phi_L(x)}{\eps}\Big) + u_{s,R}\Big(\frac{\phi_R(x)}{\eps}\Big),
\end{align}
where $\phi_L(x) = x-a$ and $\phi_R(x) = x-b$ are level-set functions associated with the left and right boundary points, and $u_{s,L}$ and $u_{s,R}$ are smooth functions constructed using standard MLPs. This representation captures the rapid transitions near the boundary points while remaining nearly constant away from them. Moreover, it follows that $u_s'' = O(1/\eps^2)$ within the boundary-layer regions.

Combining the regular and singular components, we obtain the following one-dimensional solution representation:
\begin{align}\label{Eq:1D_sol}
u(x) = u_r(x) + u_{s,L}\Big(\frac{\phi_L(x)}{\eps}\Big) + u_{s,R}\Big(\frac{\phi_R(x)}{\eps}\Big).
\end{align}
Accordingly, the neural network approximation consists of three MLPs: one MLP for the regular component $u_r$, and two MLPs for the singular components $u_{s,L}$ and $u_{s,R}$. The architectures of these MLPs will be discussed later in Subsec.~\ref{subsec:MLP}.

\paragraph{Remark}
We emphasize that the general solution ansatz in Eq.~(\ref{Eq:1D_sol}) does not require prior knowledge of the boundary layer location. Even when the boundary layer occurs only at one endpoint, the formulation remains effective and is capable of automatically identifying the boundary-layer region during training, as demonstrated by the numerical results in Sec.~\ref{sec:results}. If the layer location is known in advance, for example at the left endpoint $x=a$, the formulation may be simplified in practice by retaining only the corresponding component and setting $u_s = u_{s,L}$, although its performance is generally comparable to that of the full ansatz in Eq.~(\ref{Eq:1D_sol}).


\subsection{Two-dimensional case}\label{subsec:2D}
We next extend the solution decomposition in Eq.~(\ref{Eq:decomp}) to two spatial dimensions. 
The decomposition is discussed for different domain geometries in what follows.

\paragraph{\textbf{Regular domain}}
In most of the literature, the boundary value problem~(\ref{Eq:PDE}) in 2D space is posed on a rectangular domain $\Omega=[a,b]\times[c,d]$.
Following the same idea as in the one-dimensional case, the regular component $u_r(x,y)$ captures smooth variations over $\Omega$, while the singular component $u_s(x,y)$ is constructed in a dimension-by-dimension manner. Specifically, we define
\begin{equation}\label{Eq:2D_singular}
\begin{aligned}
u_s(x,y) = \;& \left[ u_{r,x}(x,y) + u_{s,L}\Big(\frac{\phi_L(x)}{\epsilon}\Big) + u_{s,R}\Big(\frac{\phi_R(x)}{\epsilon}\Big) \right] \\
              & \times \left[ u_{r,y}(x,y) + u_{s,B}\Big(\frac{\phi_B(y)}{\epsilon}\Big) + u_{s,T}\Big(\frac{\phi_T(y)}{\epsilon}\Big) \right],
\end{aligned}
\end{equation}
where $\phi_B(y) = y-c$ and $\phi_T(y) = y-d$ are the level-set functions corresponding to the bottom and top boundaries, with the subscripts $B$ and $T$, respectively. Here, $u_{r,x}$ and $u_{r,y}$ represent the regular components associated with the $x$- and $y$-directions, while the functions $u_{s,L}, u_{s,R}, u_{s,B}$, and $u_{s,T}$ are introduced to capture the boundary-layer behavior normal to the respective boundaries. Notably, this multiplicative construction enables $u_s$ to represent singular behavior along each coordinate direction in a separable manner and to effectively capture corner-layer phenomena \cite{RST08, E79}.

As a result, given the above construction of the singular component, the two-dimensional solution is represented as 
\begin{equation}\label{Eq:2D_sol}
\begin{aligned}
u(x,y)  = u_r(x,y) + & \left[ u_{r,x}(x,y) + u_{s,L}\Big(\frac{\phi_L(x)}{\epsilon}\Big) + u_{s,R}\Big(\frac{\phi_R(x)}{\epsilon}\Big) \right] \\
&\times
   \left[ u_{r,y}(x,y) + u_{s,B}\Big(\frac{\phi_B(y)}{\epsilon}\Big) + u_{s,T}\Big(\frac{\phi_T(y)}{\epsilon}\Big)  \right].
\end{aligned}
\end{equation}
Accordingly, multiple MLPs are employed to approximate the regular and singular components through an appropriately designed loss formulation. Analogous to the one-dimensional setting, the general two-dimensional solution representation in Eq.~\eqref{Eq:2D_sol} is capable of automatically and robustly identifying boundary-layer locations, a feature that will be further substantiated by the numerical results.

\paragraph{Remark} As in the one-dimensional case, the general formulation does not require prior knowledge of the boundary-layer locations. Boundary layers may arise only on a subset of the domain boundary rather than along the entire boundary, and the proposed ansatz remains effective in such situations. If the layer locations are known \emph{a priori} from physical considerations, the decomposition may be simplified accordingly. For example, if boundary layers are present only at $x = a$ and $x = b$, a simplified decomposition $u = u_r + u_{s,L} + u_{s,R}$ can be adopted, which is sufficient for practical implementations and in practice leads to performance comparable to that of the general ansatz in Eq.~(\ref{Eq:2D_sol}).

\paragraph{\textbf{Irregular domain}}
The neural-network-based learning framework possesses the inherent advantage of being mesh-free, allowing it to handle problems defined on complex or irregular domains where conventional numerical methods often encounter significant challenges. To leverage this flexibility within the proposed formulation, we represent the irregular domain using a level-set description. Specifically, let $\phi(x,y)$ be a smooth level-set function whose zero-level set coincides with the boundary of the irregular domain $\Omega$. Then, the two-dimensional solution on the irregular domain can be written as
\begin{align}\label{Eq:2D_sol_irregular}
u(x,y) = u_r(x,y) + u_s\Big(x,y,\frac{\phi(x,y)}{\eps}\Big),
\end{align}
where $u_r$ and $u_s$ are parameterized by standard MLPs. We emphasize that, in contrast to the singular component constructed above, which remains a univariate function of the level-set variable, the singular component $u_s$ here becomes a multivariate function of the spatial variables together with the level-set variable. This additional flexibility allows the model to automatically detect the boundary-layer regions during training for irregular domains.
From the above solution representation, let $z = \phi/\eps$, the gradient can be expressed as
\begin{align*}
\nabla u(x,y) = \nabla u_r(x,y) + \grad_\bx u_s\Big(x,y,\frac{\phi(x,y)}{\eps}\Big) +  \partial_z u_s\Big(x,y,\frac{\phi(x,y)}{\eps}\Big) \frac{\nabla \phi(x,y)}{\eps},
\end{align*}
where $\grad_\bx$ denotes the gradient with respect to the spatial variables $x$ and $y$.
By construction, the level-set function is chosen such that $\|\nabla \phi(x,y)\| = O(1)$ throughout the domain $\Omega$. As a result, the solution gradient naturally exhibits a magnitude of $O(1/\eps)$ in the vicinity of the boundary-layer regions, while remaining smooth away from the layers. Accordingly, the second-order partial derivatives scale as $O(1/\eps^2)$ within the boundary-layer regions and exhibit regular behavior elsewhere. This condition is automatically satisfied in the case of regular domains, where distance-based level-set functions are employed (see Eq.~(\ref{Eq:1D_singular}) for the one-dimensional case and Eq.~(\ref{Eq:2D_singular}) for the two-dimensional case), for which $\|\nabla \phi\| = 1$ holds identically.

\section{Implementation details}\label{sec:implementation}
In this section, we provide the implementation details of the proposed approach. Specifically, we describe the neural-network architectures employed in the numerical experiments, as well as the construction of the weighted loss function and the optimization strategy employed for model training.

\subsection{Neural network architecture}\label{subsec:MLP}

We outline the neural-network architectures designed to implement the proposed solution decomposition, with specific configurations for the one- and two-dimensional cases detailed below.

\subsubsection{One-dimensional case}
\begin{figure}[h]
\centering
\includegraphics[width=0.5\textwidth]{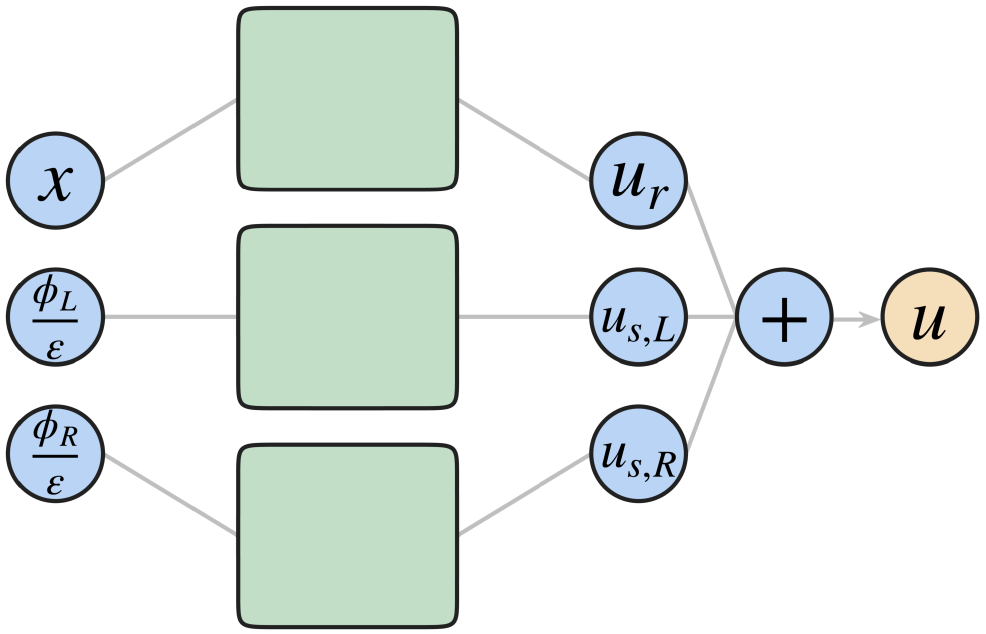}
\caption{Schematic illustration of the neural network architecture for the one-dimensional case. 
The network takes three inputs, $x$, $\phi_L/\eps$, and $\phi_R/\eps$, 
which are fed into three separate fully connected blocks with a variable number of hidden layers. 
The outputs of the hidden layers are subsequently combined and mapped to a single scalar output $u$.}
\label{Fig:architecture_1d}
\end{figure}

As discussed in Subsec.~\ref{subsec:1D}, the general solution formulation in Eq.~(\ref{Eq:1D_sol}) involves three univariate functions, namely $u_r$, $u_{s,L}$, and $u_{s,R}$, which are to be determined. A straightforward approach is to employ three separate neural networks to approximate these functions. Here, however, we adopt a more compact representation by embedding all components into a single neural network architecture, as illustrated in Fig.~\ref{Fig:architecture_1d}. Specifically, we construct the neural network whose input layer has dimension three, corresponding to the variables $x$, $\phi_L/\eps$, and $\phi_R/\eps$. These inputs are processed by three independent network blocks, each comprising a variable number of fully connected hidden layers. The features extracted from all blocks are subsequently aggregated and mapped to a single scalar value, $u$, via the final output layer.

\subsubsection{Two-dimensional case}
As discussed in Subsec.~\ref{subsec:2D}, the general representation of the two-dimensional solution admits different forms depending on the geometry of the computational domain. We next describe the corresponding neural-network architectures for the two-dimensional setting.

\begin{figure}[h]
\centering
\includegraphics[width=0.6\textwidth]{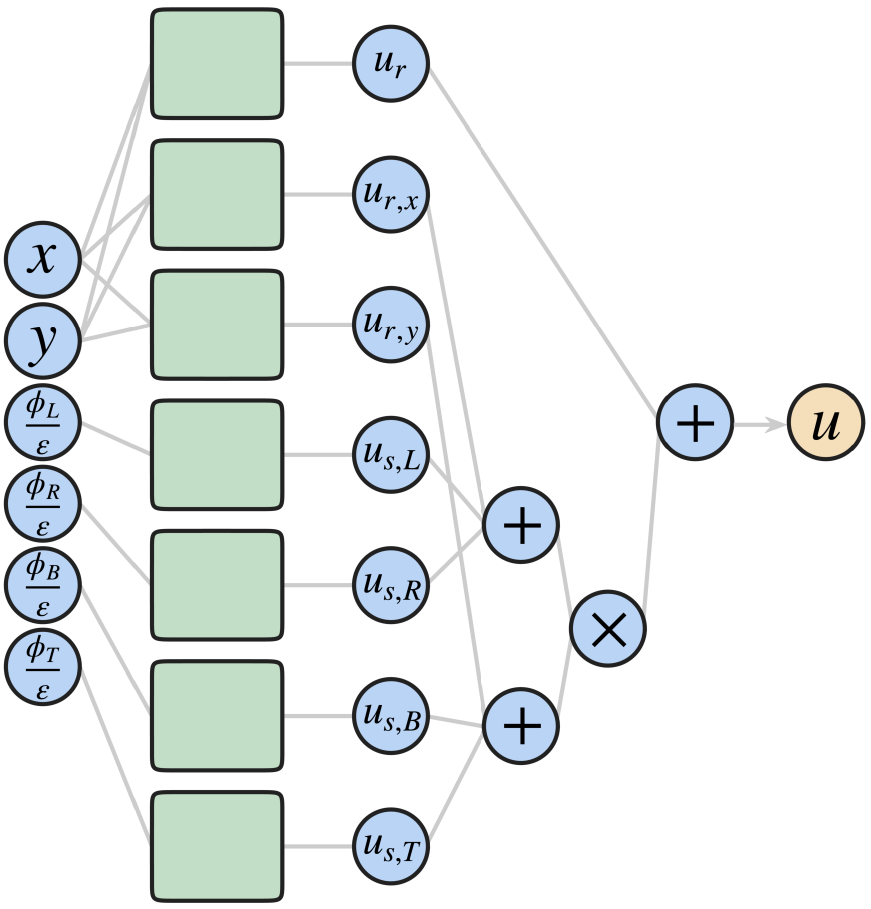}
\caption{Neural network architecture for the two-dimensional regular domain. Six inputs, $x$, $y$, and the scaled level-set functions $\phi_L/\eps$, $\phi_R/\eps$, $\phi_B/\eps$, and $\phi_T/\eps$, are processed by seven fully connected blocks to realize the solution decomposition, and the final output $u$ is obtained through algebraic combinations of the block outputs.}
\label{Fig:architecture_2d}
\end{figure}

\paragraph{\textbf{Regular domain}}
For regular domains, we consider a rectangular domain $\Omega = [a,b]\times[c,d]$ together with the solution representation in Eq.~(\ref{Eq:2D_sol}). Similar to the one-dimensional case, this decomposition is implemented via a single neural network featuring six inputs: $x$, $y$, $\phi_L/\eps$, $\phi_R/\eps$, $\phi_B/\eps$, and $\phi_T/\eps$. The network consists of seven independent fully connected blocks, each with a tunable number of hidden layers, dedicated to approximating the individual components of the decomposition. The final solution $u$ is obtained through simple algebraic combinations of the outputs from these blocks, as illustrated in Fig.~\ref{Fig:architecture_2d}.

In the one- and two-dimensional cases discussed above, the neural network is constructed to produce a single scalar output corresponding to the solution $u$. For systems of PDEs with vector-valued solutions, $\bu = [u_1, u_2, \ldots, u_n]^\top$, the same network architecture is adopted with a straightforward modification: each fully connected block is designed to output a vector of dimension $n$ instead of a scalar. The final vector-valued solution is then obtained by applying the same algebraic operations to the block outputs in a componentwise manner. This construction enables a unified neural-network architecture to be employed for both scalar and system PDEs.

\begin{figure}[h]
\centering
\includegraphics[width=0.5\textwidth]{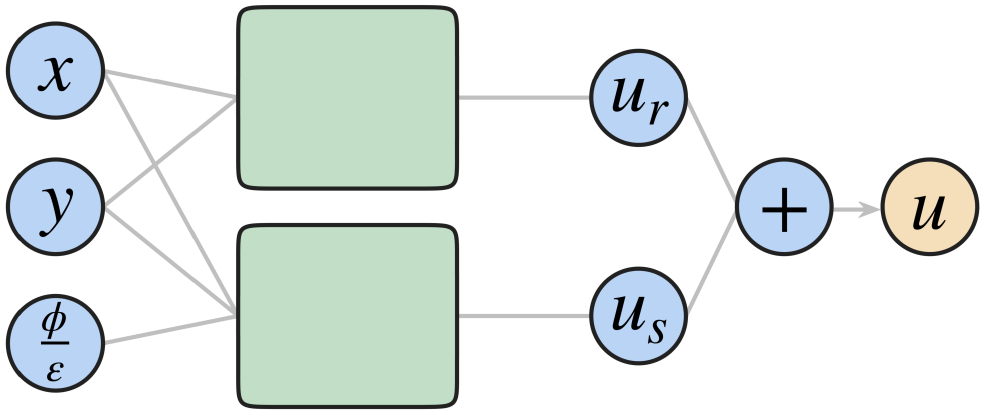}
\caption{Neural network architecture for the two-dimensional irregular domain.
The inputs consist of $x$, $y$, and the level-set function $\phi/\eps$, and two fully connected blocks are combined additively to produce the solution.}
\label{Fig:architecture_2d_ir}
\end{figure}

\paragraph{\textbf{Irregular domain}}
In the irregular-domain case, the boundary layer is assumed to occur along the boundary of the domain. Using the solution representation in Eq.~(\ref{Eq:2D_sol_irregular}), we construct a neural network architecture, with inputs $x$, $y$, and $\phi/\eps$, and two fully connected blocks whose outputs are combined additively to obtain the solution, as illustrated in Fig.~\ref{Fig:architecture_2d_ir}. 

\subsection{Weighted loss function}
We have stated the considered model problems in Subsec.~\ref{subsec:model}, all of which can be cast into the general operator form given in Eq.~(\ref{Eq:PDE}). For the system of PDEs, the governing equations admit an analogous operator form,
$\mathcal{L}_\eps(\bu) = \mathbf{f}$ in $\Omega$, with associated boundary condition $\bu = \mathbf{g}$ on $\partial\Omega$,
where $\mathbf{g} = [g_1,g_2,\cdots,g_n]^\top$ gives the prescribed boundary values.

We are now in a position to train the learnable parameters of the proposed neural network.
To this end, we employ a PINN framework~\cite{raissi2019} together with a weighted loss function tailored for singularly perturbed problems.
Let $\mathbf{x}$ denote the spatial variable in dimension $d=1,2$.
We sample a set of interior training points $\{\mathbf{x}^i\}_{i=1}^m$ in the domain $\Omega$
and a set of boundary points $\{\mathbf{x}_b^j\}_{j=1}^{m_b}$ on the domain boundary $\partial\Omega$.
The solution is obtained by minimizing a weighted mean-squared loss constructed from the residuals of the governing equation and the boundary condition in Eq.~(\ref{Eq:PDE}), namely, 
\begin{align}\label{Eq:loss}
\mathcal{J}(\boldsymbol{\theta})
=
\frac{1}{m}\sum_{i=1}^m
w(\mathbf{x}^i)
\bigl(\mathcal{L}_\eps u(\mathbf{x}^i) - f(\mathbf{x}^i)\bigr)^2
+
\frac{1}{m_b}\sum_{j=1}^{m_b}
\bigl(u(\mathbf{x}_b^j) - g(\mathbf{x}_b^j)\bigr)^2,
\end{align}
where $\boldsymbol{\theta}$ denotes the set of learnable parameters in the neural network approximation of $u$, and
$w$ is a nonnegative weight function tailored to the type of singularly perturbed PDE under consideration. All spatial derivatives appearing in the differential operator $\mathcal{L}_\eps$ are evaluated via automatic differentiation (AD)~\cite{BPRS18, GW08}.
For the system of PDEs, the loss function is defined in a componentwise manner, analogous to the scalar loss formulation described above.

%

The purpose of the weight function $w$ is to rescale all dominant terms in the underlying PDE so that their contributions are of comparable order, namely $O(1)$.
Our analysis in Section~\ref{sec:NN} shows that, with appropriately chosen level-set functions in the singular component $u_s$, the first- and second-order derivatives of the proposed neural network solution attain magnitudes of $O(1/\eps)$ and $O(1/\eps^2)$, respectively, within the boundary-layer regions, while remaining $O(1)$ away from the boundary layers.

For reaction--diffusion and Poisson--Boltzmann equations, this scaling naturally balances the dominant terms (together with their associated coefficients) in the governing equations. Consequently, no additional rescaling is required, and we simply set $w(\bx)=1$ in $\Omega$.


By contrast, for both scalar and system convection--diffusion--reaction equations, 
the dominant convection and diffusion terms scale as $O(1/\eps)$ within the 
boundary-layer regions. Consequently, the residual square in these regions scales as $O(1/\eps^2)$. To balance the contributions in the loss 
function at order $O(1)$, an additional weighting is required. 
Motivated by this observation, the weight function $w$ is designed to be of magnitude $O(\eps^2)$ in the boundary-layer regions while remaining $O(1)$ outside the boundary layers. More generally, the design of $w$ only needs to satisfy the scaling requirement discussed above. A natural choice is to construct $w$ based on the distance to the domain boundary, although the proposed framework is not restricted to this particular form.
For one-dimensional problems, we set the weight function
\begin{align*}
w(x) = \bigl(\min(x-a,\, b-x)\bigr)^2,
\end{align*}
while for two-dimensional rectangular domains, we define
\begin{align*}
w(x,y) = \bigl(\min(x-a,\, b-x,\, y-c,\, d-y)\bigr)^2.
\end{align*}
This choice of $w$ indeed reflects the characteristic boundary-layer thickness and ensures a balanced contribution of all dominant terms in the loss function.

\subsection{Sampling and training strategy}
It is essential that the solution decomposition in Eq.~(\ref{Eq:decomp})
captures the localized sharp variations within the boundary layers through
the singular component $u_s$, while the regular component $u_r$ accounts for
the smooth behavior of the solution. Both components are mainly resolved by the
PDE-residual term in the weighted loss function (the first term in
Eq.~(\ref{Eq:loss})), which is enforced both within the boundary-layer
regions and away from them. Accordingly, the training set
$\{\bx^i\}_{i=1}^m$ must consist of collocation points sampled in both the
boundary-layer regions and the interior of the domain $\Omega$.

The collocation points located near the boundary layers are primarily
responsible for identifying the singular component $u_s$ and are drawn
from truncated normal distributions with standard deviation
$\sigma_{\mbox{std}} = O(\eps)$. For example, in the one-dimensional case with
$\Omega = [a,b]$, boundary-layer points are sampled from
$\mathcal{N}(a,\sigma_{\mbox{std}}^2)$ and $\mathcal{N}(b,\sigma_{\mbox{std}}^2)$, and samples lying
outside the domain are discarded. Notably, these collocation points are also
placed near boundaries where no boundary layer is present, ensuring that
the absence of singular behavior is properly captured during training.
In contrast, the interior collocation points, which are uniformly sampled over $\Omega$, primarily contribute to resolving the regular component $u_r$. 
The boundary collocation points $\{\bx_b^j\}_{j=1}^{m_b}$ are likewise sampled uniformly along the boundary.
The same sampling strategy is also applied to the two-dimensional cases.

With the above sampling strategy, the weighted loss function defined in
Eq.~(\ref{Eq:loss}) naturally leads to a nonlinear least-squares problem,
in which the parameters of the neural networks approximating $u_r$ and
$u_s$ are determined by minimizing the aggregated PDE residuals and
boundary condition losses evaluated at the collocation points. This optimization problem can be
efficiently solved using the Levenberg--Marquardt (LM) method~\cite{Marquardt63},
which has been shown to exhibit fast and robust convergence for
PDE-based learning problems with moderate number of learnable parameters~\cite{HLL22, HLTL23, TLHL23, HSLL24}.
\subsection{Some remarks}
We conclude this section with several remarks that provide further insight
into the behavior and robustness of the proposed method.
\paragraph{\textbf{Auto detection of boundary layers}}
Although the proposed formulation does not prescribe the locations of the
boundary layers \emph{a priori}, their positions are identified automatically
during training. This capability can be attributed to the design of the
neural-network-based singular components, in which the scaled input
$\phi/\eps$ is fed into the first fully connected layer through an
affine transformation of the form $W(\phi/\eps)+b$. 
Here, $\phi$ generically refers to the level-set function associated with
the singular component under consideration, while $W$ and $b$ denote the
corresponding weights and biases. As a result, the zero level set associated with this affine
transformation in the singular component
is no longer fixed at $\{\phi=0\}$, but is instead determined implicitly by
\[
W(\phi/\eps)+b = 0,
\]
which is equivalent to $\phi(\bx)=-(b/W)\eps$ when $W\neq 0$.
Provided that the weights and biases remain $O(1)$ during training, this
relation obviously represents an $O(\eps)$ perturbation of the original zero
level set. Consequently, the training process can adaptively shift the effective
boundary-layer location within an $\eps$-neighborhood of
$\{\phi=0\}$. In regions where no boundary layer is present, the same
$\eps$-scale perturbation instead shifts the effective zero level
set outside the domain. In particular, this $\eps$-scale shift leads to no sharp solution in $\Omega$ when
no rapid transition is required by the PDE. 
Under these circumstances, activating the singular component does not
reduce the PDE residual, and the optimization process naturally drives
the singular network toward a nearly constant output over $\Omega$. The
singular component $u_s$ thus effectively degenerates to a constant in
regions without boundary layers.

\paragraph{\textbf{Rounding error effect}}
We next comment on the rounding error effects associated with the
evaluation of derivatives in singularly perturbed problems. In the proposed
formulation, all derivatives appearing in the PDE residual are computed
via AD. By analytically propagating exact
chain-rule operations through the computational graph, AD avoids the
subtraction-induced cancellation and step-size sensitivity that commonly
lead to rounding errors in traditional numerical differentiation
methods.

More importantly, within the proposed decomposition framework, the use
of appropriately scaled singular components ensures that the quantities
entering the PDE residual remain well balanced in magnitude, even for
small values of $\eps$. As a result, the derivatives produced by
AD do not involve ill-conditioned difference quotients, and the
accumulation of rounding errors is effectively mitigated during training.

\section{Numerical results}\label{sec:results}
In this section, we present a series of numerical experiments to assess
the robustness and accuracy of the proposed method. Throughout all
examples, a single-hidden-layer architecture is employed in each network
block with sigmoid activation functions, resulting in a moderate number
of trainable parameters ranging from several hundred to approximately one
thousand. The LM optimizer is applied until one of
the stopping criteria is met, namely, the loss $\mathcal{J}<10^{-15}$ or a
maximum of 2000 iterations.

To quantify the solution accuracy, we report the mean relative $L^2$ and
$L^\infty$ errors computed over up to five independent trial runs, where
$u_{\mathcal N}$ denotes the neural-network approximation of the solution
and $u$ denotes the corresponding exact solution. The errors $\|u-u_{\mathcal N}\|/\|u\|$ are
evaluated at test points
located in the interior of the domain as well as in the boundary-layer
regions. These test points are sampled using the same strategy as that
adopted for the training collocation points, with the total number being twice that of the training points.

\subsection{One-dimensional problems}
We begin with a set of one-dimensional examples in the domain $\Omega=[0,1]$. In all one-dimensional experiments, we employ 50 neurons in each network
block for both scalar- and vector-valued solutions.
This results in 450 trainable parameters for Examples 1–2 and 600 trainable parameters for Example 3, respectively.
A total of $m=1500$ collocation points are used to evaluate the
PDE residual, consisting of 500 points uniformly sampled in the interior
of the domain and 1000 points distributed near the two endpoints (500 near
each endpoint) to capture possible boundary-layer behavior. In addition,
two boundary points ($m_b=2$) are imposed separately to enforce the
boundary conditions.

\paragraph{\textbf{Example~1}}
We consider the one-dimensional convection--diffusion--reaction problem:
\begin{align*}
-\eps u''(x) - (1+\eps)u'(x) - u(x) = 0,
\end{align*}
subject to the boundary conditions $u(0)=0$ and $u(1)=1$. The analytic solution of this problem is given by
\begin{align*}
u(x) = \frac{\exp(-x)-\exp(-x/\eps)}{\exp(-1)-\exp(-1/\eps)},
\end{align*}
from which it is evident that a single boundary layer is present at the left endpoint $x=0$.

\begin{figure}[h]
\centering
\includegraphics[width=0.9\textwidth]{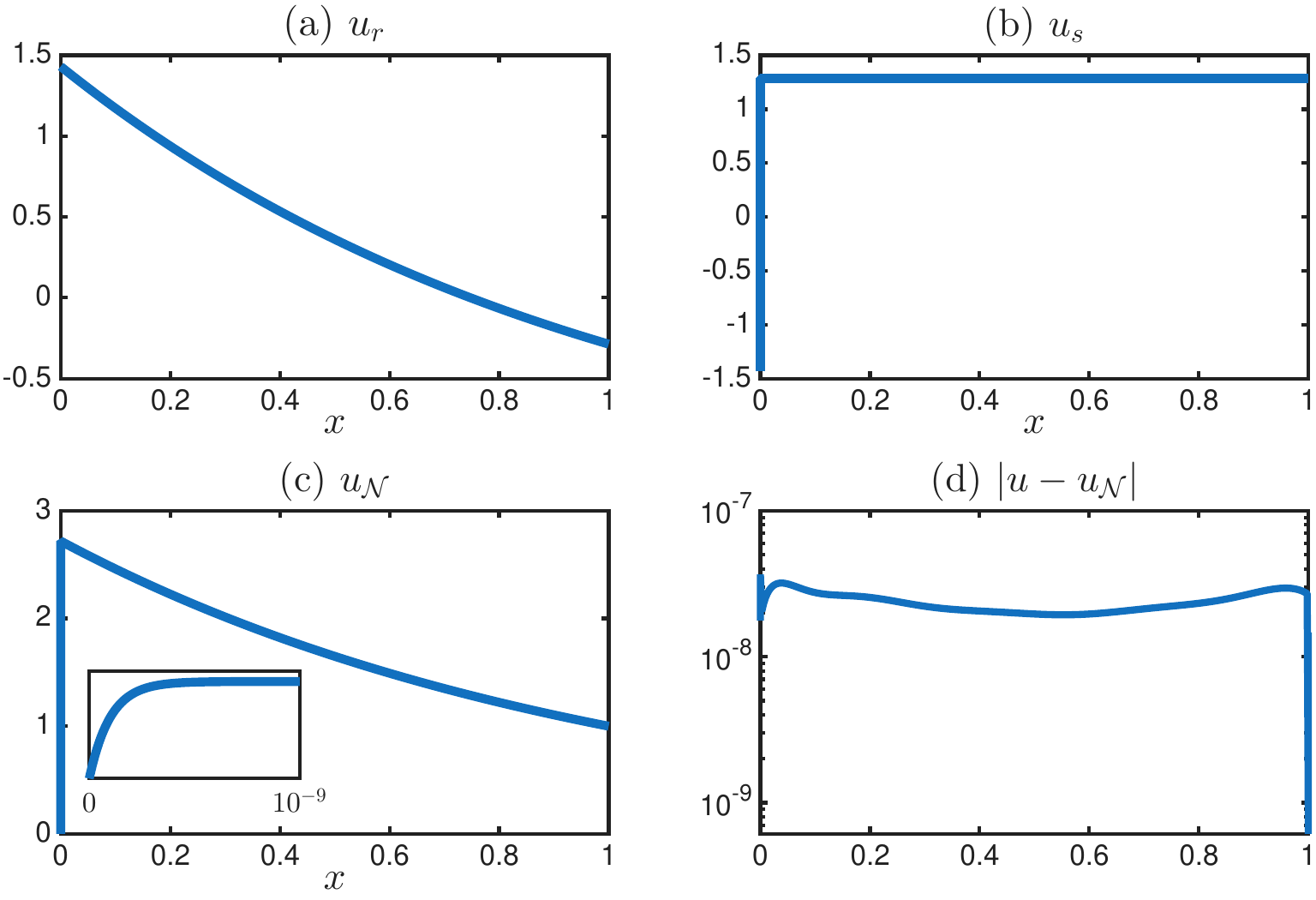}
\caption{One-dimensional convection–diffusion–reaction problem with
$\eps=10^{-10}$ in Example 1.
(a) Learned regular component $u_r$, which remains smooth over the entire
domain.
(b) Learned singular component $u_s$, accurately capturing the boundary
layer near the left endpoint while remaining nearly constant elsewhere.
(c) Neural-network approximation $u_{\mathcal N}$ with an zoom-in inset view near the boundary layer at $x=0$.
(d) Absolute error $|u-u_{\mathcal N}|$, showing that the maximum error is
localized near the boundary layer and remains of comparable magnitude
throughout the domain.}
\label{Fig:EX1}
\end{figure}

We perform a numerical experiment with $\eps=10^{-10}$. In Fig.~\ref{Fig:EX1}(a--b), we display the learned regular and singular components, respectively. As can be seen, the regular component remains smooth over the entire domain, whereas the singular component accurately captures the boundary-layer structure near the left endpoint. On the right side of the domain, where no boundary layer is present, the singular component remains nearly constant. The neural-network approximation $u_{\mathcal N}$ is shown in Fig.~\ref{Fig:EX1}(c), where the inset provides a zoomed-in view of the boundary layer. The corresponding absolute error $|u-u_{\mathcal N}|$ is reported in Fig.~\ref{Fig:EX1}(d). It can be observed that the error magnitude remains on the order of $10^{-8}$ even for such a small value of $\eps$, a level of accuracy that is challenging to achieve using traditional numerical methods. Although the maximum error occurs near the boundary layer at $x=0$, its magnitude remains comparable to the error level observed away from the boundary-layer region.

Moreover, Table~\ref{Table:EX1} presents a systematic investigation of the
proposed method for various values of $\eps$. It is observed that
the method maintains a consistently high level of accuracy as
$\eps$ decreases from $10^{-2}$ to $10^{-10}$. In particular,
neither the relative $L^2$ nor the $L^\infty$ error exhibits any
noticeable growth in the singularly perturbed regime. Instead, both error
measures remain approximately $10^{-7}$–$10^{-8}$, indicating
that the presence of an increasingly thin boundary layer does not degrade
the numerical accuracy.

In all cases, the training process reduces the loss value to
approximately $\mathcal{J}(\boldsymbol{\theta}) \approx 10^{-15}$.
Although a rigorous theoretical connection between the loss value and
the solution error is not established here, the observed errors
consistently remain within the scale of $\sqrt{\mathcal{J}}$, suggesting
a strong correlation between the minimization of the loss function and
the achieved solution accuracy.

\begin{table}[!h]
\centering
\begin{tabular}{c|ccccc}
\hline
$\eps$ & $10^{-2}$ & $10^{-4}$ & $10^{-6}$ & $10^{-8}$ & $10^{-10}$ \\
\hline
rel. $L^2$ error        & 2.43E$-$08 & 8.69E$-$09 & 4.46E$-$09 & 1.06E$-$08 & 1.38E$-$08 \\
rel. $L^\infty$ error  & 2.45E$-$08 & 2.91E$-$08 & 1.59E$-$08 & 1.57E$-$08 & 1.80E$-$08 \\
\hline
\end{tabular}
\caption{Relative $L^2$ and $L^\infty$ errors for the one-dimensional
convection–diffusion–reaction problem in Example~1 with varying values of
$\eps$. The results demonstrate that the proposed method maintains
a consistently high level of accuracy as $\eps$ decreases from
$10^{-2}$ to $10^{-10}$, with no noticeable error growth in the
singularly perturbed regime.}
\label{Table:EX1}
\end{table}

\paragraph{\textbf{Example~2}}
In the second example, we examine a one-dimensional reaction--diffusion
problem given by
\begin{align*}
-\eps^{2} u''(x) + u(x) = 0,
\end{align*}
with the two-endpoint boundary conditions $u(0)=0$ and $u(1)=1$. The analytic solution is
\begin{align*}
u(x) = \frac{-\exp\!\big(-(1-x)/\eps\big)
+ \exp\!\big(-(1+x)/\eps\big)}{\exp(-2/\eps)-1},
\end{align*}
which exhibits a boundary layer at the right endpoint $x=1$. 

\begin{figure}[t]
\centering
\includegraphics[width=0.9\textwidth]{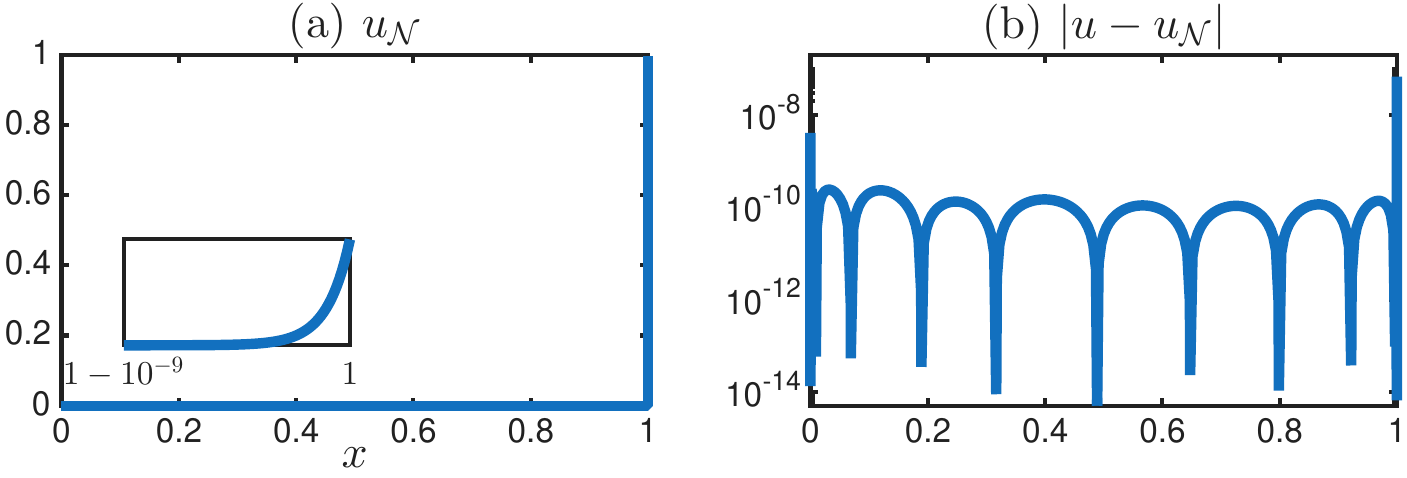}
\caption{One-dimensional reaction–diffusion problem in Example~2 with
$\eps=10^{-10}$.
(a) Neural-network approximation $u_{\mathcal N}$. An inset provides a zoomed-in view near $x=1$ to highlight the
boundary-layer structure.
(b) Absolute error $|u-u_{\mathcal N}|$.}
\label{Fig:EX2}
\end{figure}

We conduct a numerical experiment with $\eps=10^{-10}$ and the
resulting neural-network approximation and its absolute error are shown
in Fig.~\ref{Fig:EX2}. It can be observed that the proposed method successfully identifies the location of the boundary layer at the right
endpoint $x=1$. Although the intrinsic physical nature of the underlying PDE differs from that of the previous convection–diffusion–reaction case,
the same network architecture and loss formulation in
Eq.~(\ref{Eq:loss}) are directly applied without modification.

The relative $L^2$ and $L^\infty$ errors for Example~2 with different
values of $\eps$ are summarized in Table~\ref{Table:EX2}. The
results show that the proposed method achieves a high-accuracy level, with error magnitudes on the order of $10^{-7}$–$10^{-8}$,
across all tested values of $\eps$ from $10^{-2}$ to $10^{-10}$. This evidence
indicates that, despite the different reaction--diffusion structure of
the underlying PDE, the proposed framework preserves its robustness with
respect to increasingly thin boundary layers.

\begin{table}[h]
\centering
\begin{tabular}{c|ccccc}
\hline
$\eps$ & $10^{-2}$ & $10^{-4}$ & $10^{-6}$ & $10^{-8}$ & $10^{-10}$ \\
\hline
rel. $L^2$ error        & 1.73E$-$07 & 2.17E$-$07 & 1.08E$-$07 & 7.42E$-$08 & 6.50E$-$08 \\
rel. $L^\infty$ error  & 5.09E$-$08 & 9.00E$-$08 & 2.74E$-$08 & 3.38E$-$08 & 3.21E$-$08 \\
\hline
\end{tabular}
\caption{Relative $L^2$ and $L^\infty$ errors for the
one-dimensional reaction--diffusion problem in Example~2 with varying
values of $\eps$. The results demonstrate that the proposed method
maintains a stable level of accuracy across a wide range of
$\eps$, including the strongly singularly perturbed regime.}
\label{Table:EX2}
\end{table}

\paragraph{\textbf{Example~3}}
Now, we turn to solve a benchmark from~\cite{HYY18} involving a coupled system of one-dimensional convection--diffusion equations as follows:
\begin{align*}
-\eps
\begin{bmatrix}
u_1'' \\
u_2''
\end{bmatrix}
-
\begin{bmatrix}
3 & -1 \\
4 & -1
\end{bmatrix}
\begin{bmatrix}
u_1' \\
u_2'
\end{bmatrix}
=
\begin{bmatrix}
2x+2 \\
4
\end{bmatrix},
\end{align*}
with the homogeneous boundary conditions $u_1(0)=u_1(1)=u_2(0)=u_2(1)=0$ such that the analytic solution reads
\begin{align*}
u_1(x) = 2W_1(x) + W_2(x),\quad u_2(x) = 4W_1(x),
\end{align*}
where
\begin{align*}
W_1(x) =\;&
\frac{e^{-1/\eps}\Big[\eps(-2x^2+2x+2)+\eps^2(8x-4)-1\Big]}
{\eps(1-e^{-1/\eps})^2}
+ \frac{e^{-x/\eps}(-2\eps x + x - 4\eps^2)}
{\eps(1-e^{-1/\eps})^2} \\[0.4em]
&+ \frac{e^{-(x+1)/\eps}\Big[(2\eps-1)x+4\eps^2-2\eps+1\Big]}
{\eps(1-e^{-1/\eps})^2}
+ \frac{e^{-2/\eps}\eps x(x-4\eps-1)
-\eps(x-1)(4\eps-x)}
{\eps(1-e^{-1/\eps})^2},
\\[0.8em]
W_2(x) =\;&
\frac{x(x-2\eps)e^{-1/\eps}
+(x-1)(-x+2\eps-1)
+(2\eps-1)e^{-x/\eps}}
{1-e^{-1/\eps}}.
\end{align*}

Both solution components $u_1$ and $u_2$ display boundary layers at the left endpoint $x=0$.
The coupled system is treated within the same training framework as the
scalar problems discussed earlier, where the loss function is constructed
by summing the weighted PDE residuals and boundary-condition residuals over all
solution components.

\begin{figure}[h]
\centering
\includegraphics[width=0.9\textwidth]{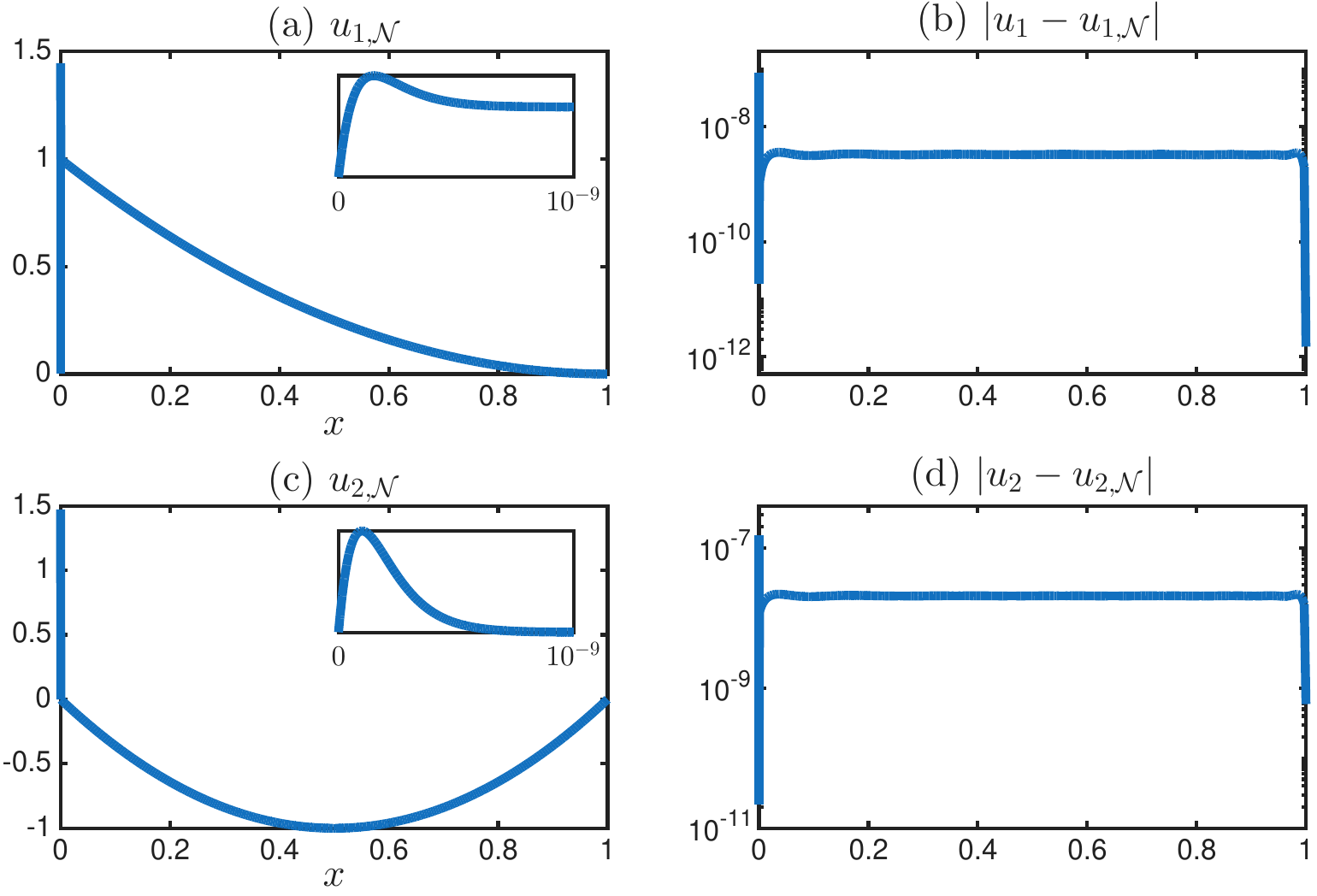}
\caption{Numerical results for the coupled system problem in Example~3 with
$\eps=10^{-10}$.
Panels (a) and (b) show the neural-network approximation $u_{1,\mathcal N}$
and its corresponding absolute error, respectively, while panels (c) and
(d) present the corresponding results for $u_{2,\mathcal N}$. Both
solution components exhibit non-monotonic boundary-layer structures near
the left endpoint $x=0$ (see the insets in panels (a) and (c) with the zoomed-in views near the left endpoint
$x=0$), which are accurately resolved by the proposed
method. The absolute errors remain at the level of $10^{-7}$--$10^{-8}$
throughout the domain.}
\label{Fig:EX3}
\end{figure}

We perform a numerical experiment with $\eps=10^{-10}$, and the results
are presented in Fig.~\ref{Fig:EX3}. Panels (a) and (b) display the
neural-network approximation $u_{1,\mathcal N}$ and its corresponding
absolute error, respectively, while panels (c) and (d) show the
corresponding results for $u_{2,\mathcal N}$. Notably, in this strongly coupled system, the boundary-layer
profiles of both solution components exhibit a non-monotonic structure,
characterized by rapid increasing and decreasing variations within the
boundary-layer region near $x=0$. These localized sharp variations are
accurately captured by the proposed method, with the magnitude of the
absolute errors remaining at the level of $10^{-7}$--$10^{-8}$.

The relative $L^2$ and $L^\infty$ errors for the one-dimensional coupled
system in Example~3 are summarized in Table~\ref{Table:EX3}. The results
indicate that both solution components are approximated with comparable
accuracy across all tested values of $\eps$. Although the two
components exhibit different boundary-layer profiles and non-monotonic
structures, the error levels for $u_1$ and $u_2$ remain of similar
magnitudes. In particular, no pronounced disparity is observed between
the $L^2$ and $L^\infty$ errors of the two components, suggesting that the
proposed framework resolves each equation in the coupled system in a
balanced manner.

\begin{table}[h]
\centering
\begin{tabular}{c|ccccc}
\hline
$\eps$ & $10^{-2}$ & $10^{-4}$ & $10^{-6}$ & $10^{-8}$ & $10^{-10}$ \\
\hline
rel. $L^2$ error $(u_1)$        & 2.12E$-$08 & 2.53E$-$08 & 3.61E$-$08 & 2.91E$-$08 & 2.06E$-$08 \\
rel. $L^\infty$ error $(u_1)$  & 1.71E$-$08 & 1.68E$-$07 & 3.94E$-$08 & 1.43E$-$07 & 2.94E$-$08 \\
rel. $L^2$ error $(u_2)$        & 3.78E$-$08 & 7.27E$-$08 & 1.69E$-$07 & 1.27E$-$07 & 7.18E$-$08 \\
rel. $L^\infty$ error $(u_2)$  & 2.51E$-$08 & 3.08E$-$07 & 8.14E$-$08 & 2.54E$-$07 & 5.99E$-$08 \\
\hline
\end{tabular}
\caption{Relative $L^2$ and $L^\infty$ errors for each
solution component of the one-dimensional coupled system in Example~3
with varying values of $\eps$.}
\label{Table:EX3}
\end{table}


\subsection{Two-dimensional problems}
We next present a series of numerical experiments for two-dimensional singularly perturbed problems to assess the performance and robustness of the proposed method. Examples~4--6 are posed on the regular domain
$\Omega = [0,1]\times[0,1]$, whereas Example~7 considers an irregular domain.

In all cases, a total of $m = 2500$ collocation points are employed, including 500 points uniformly distributed in the interior and 2000 points placed in the boundary-layer regions. For the regular domain, the boundary-layer points are evenly distributed along the four sides, with 500 points allocated to each side.
In addition, $m_b = 880$ ($\approx 40\sqrt{m}$) boundary points are uniformly sampled for the regular domain, while $m_b = 220$ ($\approx 10\sqrt{m}$) boundary points are used in the irregular-domain case. We employ one-hidden layer structure with 50 neurons in each network block for scalar problems in Example~4--5, and 35 neurons for the system problem considered in Example~6.

\paragraph{\textbf{Example 4}}
The first two-dimensional test is taken from~\cite{HY16}, where we consider a
convection--diffusion--reaction equation of the form
\begin{align*}
-\eps\laplace u + \mathbf{a}\cdot\grad u+ u = f,
\end{align*}
subject to homogeneous Dirichlet boundary conditions
$u=0$ on $\partial\Omega$. We choose a constant velocity field $\mathbf{a} = [a_1,a_2]^\top = [1/2,\sqrt{3}/2]^\top$, for which an
analytic solution is available and given by
\begin{align*}
u(x,y)
&=
\left[
\frac{x^2}{2a_1}
+ \frac{\eps x}{a_1^2}
+ \left(\frac{1}{2a_1} + \frac{\eps}{a_1^2}\right)
\frac{e^{-a_1/\eps} - e^{-a_1(1-x)/\eps}}{1 - e^{-a_1/\eps}}
\right]
\\
&\quad\times
\left[
\frac{y^2}{2a_2}
+ \frac{\eps y}{a_2^2}
+ \left(\frac{1}{2a_2} + \frac{\eps}{a_2^2}\right)
\frac{e^{-a_2/\eps} - e^{-a_2(1-y)/\eps}}{1 - e^{-a_2/\eps}}
\right].
\end{align*}

The corresponding source term $f$ is obtained by substituting the analytic solution into the governing equation. Since the prescribed velocity field is directed along the positive $x$- and $y$-directions, the solution exhibits boundary-layer behavior along the outflow sides $x = 1$ and $y = 1$.

\begin{figure}[h]
\centering
\includegraphics[width=0.9\textwidth]{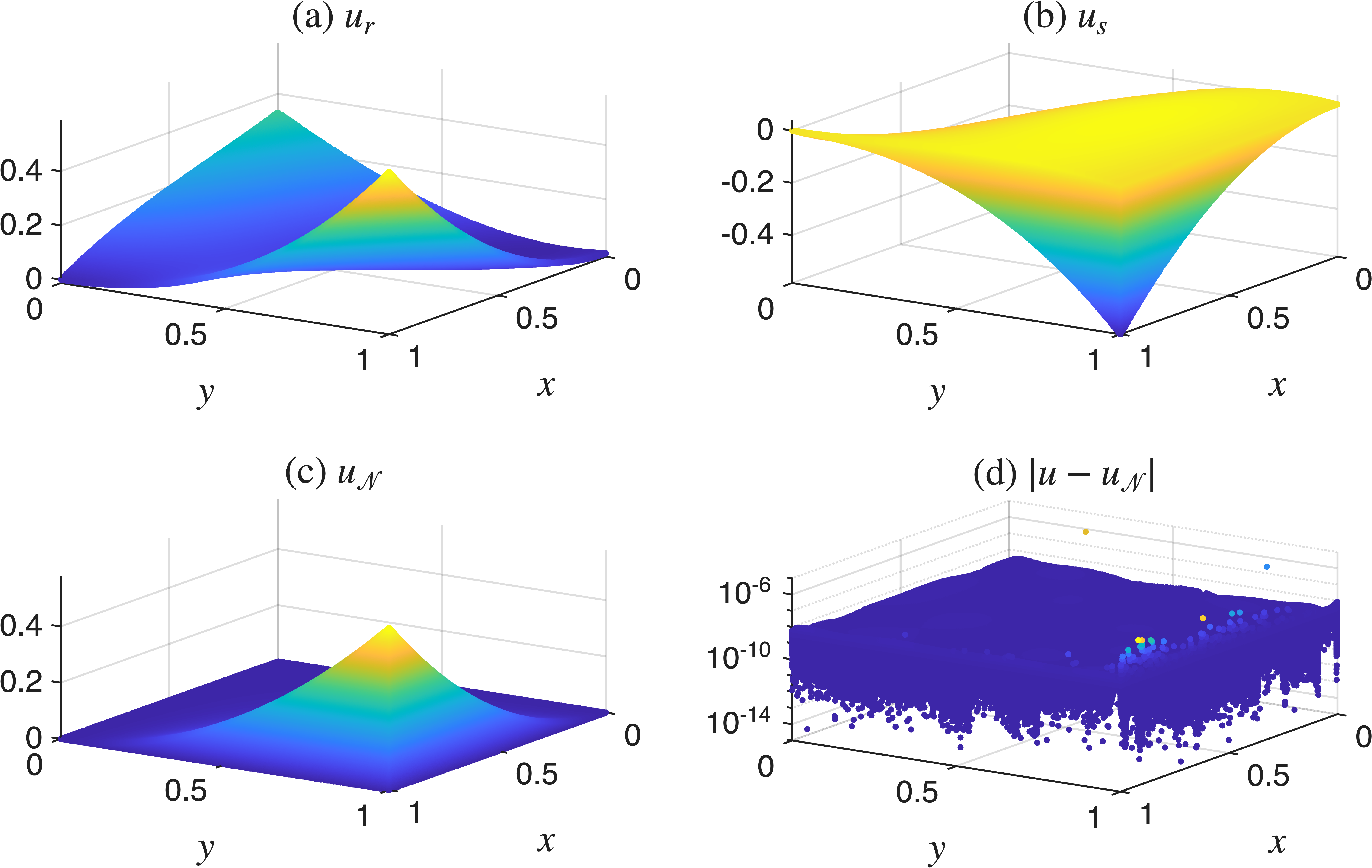}
\caption{Numerical results for Example~4 with $\eps = 10^{-10}$. Panel (a) shows the regular component $u_r$, which remains smooth over the entire domain. Panel (b) displays the singular component $u_s$, exhibiting sharp boundary-layer behavior along $x=1$ and $y=1$. Panel (c) presents the predicted solution $u_\mathcal{N}$, and panel (d) shows the corresponding absolute error.}
\label{Fig:EX4}
\end{figure}

With the adopted neural network architecture, the total number of learnable parameters is 1200. We now present the numerical results for $\eps = 10^{-10}$ in Fig.~\ref{Fig:EX4}. As shown in panels (a) and (b), the regular component $u_r$ remains smooth over the entire domain, whereas the singular component $u_s$ displays sharp localized variations along the boundary-layer regions nearby $x=1$ and $y=1$, showing that our method successfully detect the location of boundary layers under the general representation solution in Eq.~(\ref{Eq:2D_sol}). Consequently, the predicted solution $u_\mathcal{N}$ in panel (c) achieves an absolute error of order $10^{-7}$, as illustrated in panel (d). Although the exact solution in this example can be written as the product of two univariate functions, consistent with the proposed general solution ansatz in Eq.~(\ref{Eq:2D_sol}), the numerical solution obtained by our method and shown in Fig.~\ref{Fig:EX4} does not reduce to a purely multiplicative form. This demonstrates that the proposed formulation is sufficiently flexible to capture solution structures beyond separable cases. 

Table~\ref{Table:EX4} reports the relative $L^2$ and $L^\infty$ errors for this example under different values of the perturbation parameter $\eps$. It can be seen that the error levels remain nearly unchanged across several orders of magnitude of $\eps$, suggesting that the dominant error source arises from the neural network approximation and training procedure rather than from the singular perturbation scale itself.

\begin{table}[h]
\centering
\begin{tabular}{c|ccccc}
\hline
$\eps$ & $10^{-2}$ & $10^{-4}$ & $10^{-6}$ & $10^{-8}$ & $10^{-10}$ \\
\hline
rel. $L^2$ error        & 4.46E$-$07 & 4.92E$-$07 & 5.05E$-$07 & 7.09E$-$07 & 4.53E$-$07 \\
rel. $L^\infty$ error  & 3.55E$-$07 & 1.89E$-$06 & 7.23E$-$07 & 8.37E$-$07 & 5.38E$-$07 \\
\hline
\end{tabular}
\caption{Relative $L^2$ and $L^\infty$ errors for the two-dimensional convection--diffusion--reaction equation in Example~4 with varying $\eps$. In each case, there are in total 1200 learnable parameters to be trained.}
\label{Table:EX4}
\end{table}

\paragraph{\textbf{Example 5}}
The second two-dimensional test is from~\cite{TD07}, where we consider a variable-coefficient convection--diffusion equation:
\begin{align*}
-\eps\laplace u + \frac{1}{1+y}u_y = \left( \frac{1}{1+y} - 2\eps \right)\exp(y-x).
\end{align*}
The source term and boundary conditions are chosen so that the analytic solution is given by
\begin{align*}
u(x,y) = \exp(y-x) + (1+y)\exp\left( \frac{1}{\eps}\ln\frac{1+y}{2} \right).
\end{align*}
Since the convection field has no horizontal component and its vertical component is strictly positive throughout the domain, the solution develops a boundary layer along the outflow boundary $y=1$.  

\begin{figure}[h]
\centering
\includegraphics[width=0.9\textwidth]{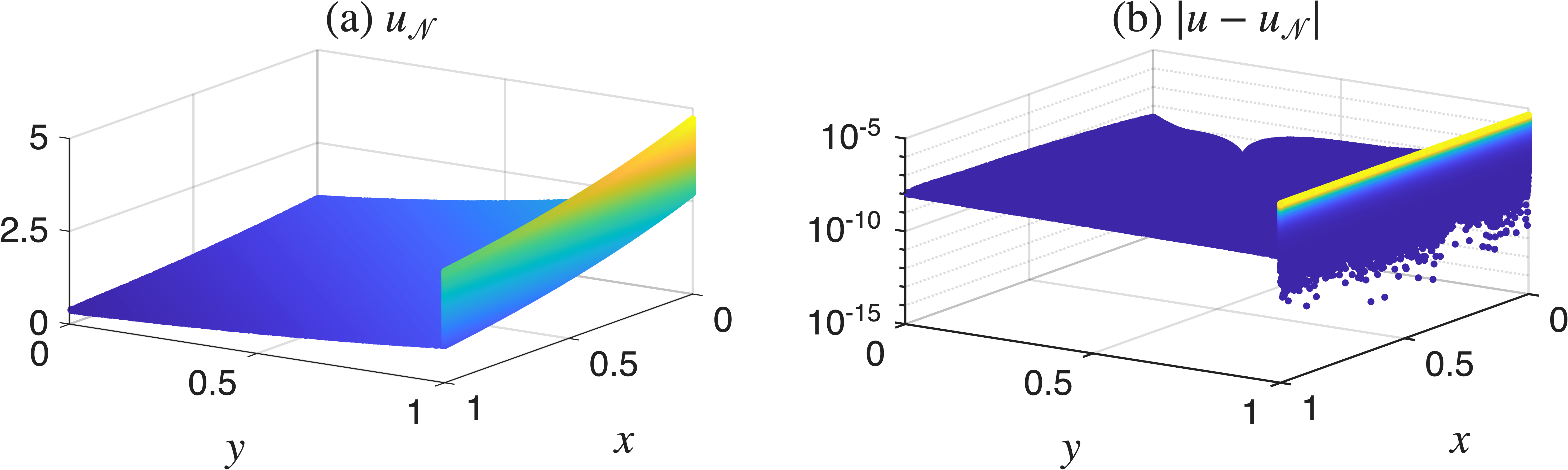}
\caption{Two-dimensional reaction–diffusion problem in Example~4 with
$\eps=10^{-10}$.
(a) Neural-network approximation $u_{\mathcal N}$. (b) Absolute error $|u-u_{\mathcal N}|$.}
\label{Fig:EX5}
\end{figure}

We solve this singularly perturbed problem using the same neural network architecture as in the previous example, with a total of 1200 learnable parameters. The neural-network prediction and the corresponding absolute error are depicted in Fig.~\ref{Fig:EX5}. These results demonstrate that the proposed method remains effective for variable-coefficient convection fields and achieves accurate predictions without any further modification of the loss formulation.

The numerical accuracy for Example~5 is summarized in Table~\ref{Table:EX5}, where the relative $L^2$ and $L^\infty$ errors are reported for different values of the perturbation parameter $\varepsilon$. In contrast to the previous example with constant convection coefficients, the present problem involves a spatially varying velocity field. Nevertheless, the error levels remain well controlled across all tested values of $\varepsilon$, indicating that the proposed method maintains its accuracy in the presence of variable-coefficient convection.

\begin{table}[h]
\centering
\begin{tabular}{c|ccccc}
\hline
$\eps$ & $10^{-2}$ & $10^{-4}$ & $10^{-6}$ & $10^{-8}$ & $10^{-10}$ \\
\hline
rel. $L^2$ error        & 7.83E$-$07 & 9.35E$-$07 & 2.17E$-$07 & 1.30E$-$07 & 1.97E$-$07 \\
rel. $L^\infty$ error  & 8.21E$-$07 & 9.95E$-$07 & 9.17E$-$07 & 5.14E$-$07 & 7.28E$-$07 \\
\hline
\end{tabular}
\caption{Relative $L^2$ and $L^\infty$ errors for Example~5 with a variable--coefficient convection field, evaluated for different values of the perturbation parameter $\varepsilon$.}
\label{Table:EX5}
\end{table}

\paragraph{\textbf{Example 6}}
Next, we extend the proposed method to solve the system of two-dimensional convection--diffusion equations referred from~\cite{HYY18}:
\begin{align*}
-\eps \Delta
\begin{bmatrix}
u_1 \\[2pt]
u_2
\end{bmatrix}
-
\begin{bmatrix}
1/2 & 1 \\[2pt]
0 & 1/2
\end{bmatrix}
\begin{bmatrix}
u_1 \\[2pt]
u_2
\end{bmatrix}_x
-
\begin{bmatrix}
\sqrt{3}/2 & 1 \\[2pt]
0 & \sqrt{3}/2
\end{bmatrix}
\begin{bmatrix}
u_1 \\[2pt]
u_2
\end{bmatrix}_y
=
\begin{bmatrix}
f_1 \\[2pt]
f_2
\end{bmatrix}.
\end{align*}
The source term and the imposed boundary conditions are given through the analytic solution:
\begin{align*}
u_1(x,y)
&=
\frac{
x^2 e^{-1/(2\eps)}
- (x^2-1)-e^{-x/(2\eps)}
+ 4\eps\!\left(x-1-xe^{-1/(2\eps)}+e^{-x/(2\eps)}\right)
}{
1-e^{-1/(2\eps)}
}
\\
&\quad\times
\frac{
\sqrt{3}\!\left(
y^2 e^{-\sqrt{3}/(2\eps)}
- (y^2-1)-e^{-\sqrt{3}y/(2\eps)}
\right)
+ 4\eps\!\left(
y-1-ye^{-\sqrt{3}/(2\eps)}+e^{-\sqrt{3}y/(2\eps)}
\right)
}{
3\!\left(1-e^{-\sqrt{3}/(2\eps)}\right)
},
\\[1ex]
u_2(x,y)
&=
-\left[
\frac{
2\!\left(
-x^3 e^{-1/(2\eps)} + x^3 - 1 + e^{-x/(2\eps)}
\right)
- 12\eps\!\left(
-x^2 e^{-1/(2\eps)} + x^2 - 1 + e^{-x/(2\eps)}
\right)
}{
3\!\left(1-e^{-1/(2\eps)}\right)
}
\right.
\\
&\qquad\left.
+
\frac{
48\eps^2\!\left(
x-1-xe^{-1/(2\eps)}+e^{-x/(2\eps)}
\right)
}{
3\!\left(1-e^{-1/(2\eps)}\right)
}
\right]
\nonumber\\
&\quad\times
\left[
\frac{
2\!\left(
-y^3 e^{-\sqrt{3}/(2\eps)} + y^3 - 1 + e^{-\sqrt{3}y/(2\eps)}
\right)
- 4\sqrt{3}\eps\!\left(
-y^2 e^{-\sqrt{3}/(2\eps)} + y^2 - 1 + e^{-\sqrt{3}y/(2\eps)}
\right)
}{
3\sqrt{3}\!\left(1-e^{-\sqrt{3}/(2\eps)}\right)
}
\right.
\nonumber\\
&\qquad\left.
+
\frac{
16\eps^2\!\left(
y-1-ye^{-\sqrt{3}/(2\eps)}+e^{-\sqrt{3}y/(2\eps)}
\right)
}{
3\sqrt{3}\!\left(1-e^{-\sqrt{3}/(2\eps)}\right)
}
\right].
\end{align*}
Both solution components display bound layers near the $x$- and $y$-axes for small values of $\eps$. 

For systems of singularly perturbed PDEs, the locations of boundary layers are in general not known \emph{a priori}, due to the coupling between different solution components and the resulting complex layer structures. This difficulty is further exacerbated when the associated coefficient matrices in a system cannot be diagonalized, as the fast and slow modes cannot be cleanly separated and no clear inflow--outflow characterization is available. The system considered in this work falls precisely into this category, making the \emph{a priori} identification of boundary-layer locations particularly challenging.

\begin{figure}[h]
\centering
\includegraphics[width=0.9\textwidth]{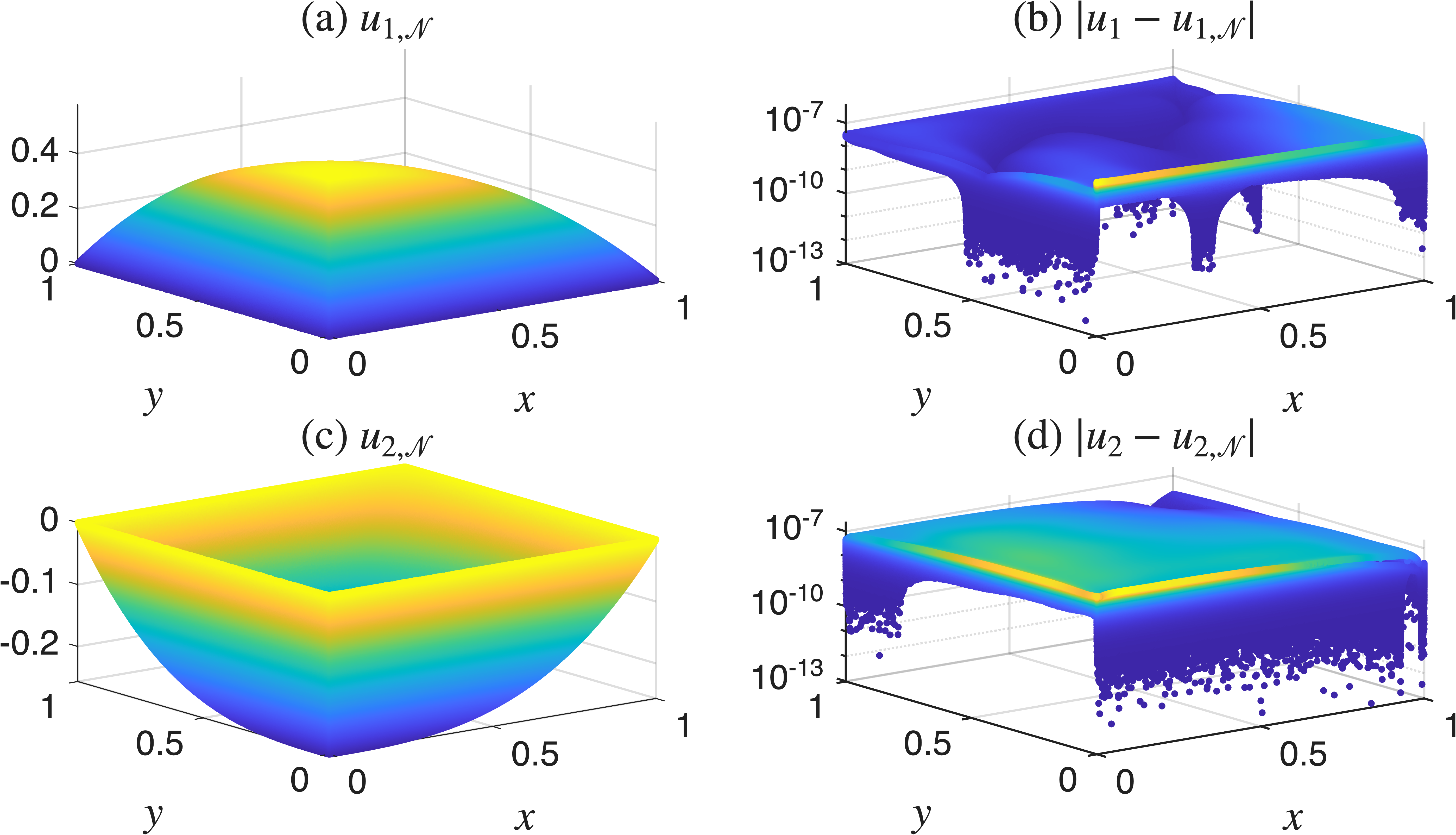}
\caption{Numerical results for the coupled system problem in Example~6 with
$\eps=10^{-10}$.
Panels (a) and (b) show the neural-network approximation $u_{1,\mathcal N}$
and its corresponding absolute error, respectively, while panels (c) and
(d) present the corresponding results for $u_{2,\mathcal N}$.}
\label{Fig:EX6}
\end{figure}

A numerical experiment is conducted with $\varepsilon = 10^{-10}$ using a neural network with 1085 learnable parameters; the neural-network predicted solution pair $[u_{1,\mathcal{N}}, u_{2,\mathcal{N}}]$ and their corresponding absolute errors are presented in Fig.~\ref{Fig:EX6}. Despite the absence of any \emph{a priori} information on the boundary-layer locations for this non-diagonalizable system, the proposed method is still able to correctly identify the boundary-layer regions for both solution components. Moreover, accurate approximations are obtained for each component, demonstrating the effectiveness of the proposed framework for singularly perturbed systems with coupled and non-diagonalizable structures.

Additional quantitative results for the system problem are reported in Table~\ref{Table:EX6}, where the relative $L^2$ and $L^\infty$ errors are shown separately for the two solution components $u_1$ and $u_2$. Although the two components are coupled and may exhibit different layer structures, both error measures remain at comparable magnitudes across all tested values of $\varepsilon$. This indicates that the proposed method provides balanced accuracy for each component of the system and does not favor one variable over the other.

\begin{table}[!h]
\centering
\begin{tabular}{c|ccccc}
\hline
$\eps$ & $10^{-2}$ & $10^{-4}$ & $10^{-6}$ & $10^{-8}$ & $10^{-10}$ \\
\hline
rel. $L^2$ error $(u_1)$        & 2.06E$-$07 & 9.49E$-$08 & 8.99E$-$07 & 6.38E$-$07 & 2.05E$-$07 \\
rel. $L^\infty$ error $(u_1)$  & 2.14E$-$07 & 8.40E$-$07 & 2.16E$-$06 & 3.03E$-$06 & 5.68E$-$07 \\
rel. $L^2$ error $(u_2)$        & 1.60E$-$07 & 1.62E$-$07 & 1.64E$-$06 & 7.87E$-$07 & 3.55E$-$07 \\
rel. $L^\infty$ error $(u_2)$  & 1.86E$-$07 & 6.53E$-$07 & 8.54E$-$06 & 2.56E$-$06 & 6.85E$-$07 \\
\hline
\end{tabular}
\caption{Relative $L^2$ and $L^\infty$ errors for the two solution components $u_1$ and $u_2$ in Example~6, evaluated for different values of the perturbation parameter $\varepsilon$.}
\label{Table:EX6}
\end{table}

\paragraph{\textbf{Example 7}}
We consider the two-dimensional Poisson--Boltzmann equation introduced in Sec.~\ref{subsec:model}, subject to the boundary condition $u = 1$ imposed along a curved irregular domain, whose zero-level-set function is given by
\begin{align*}
\phi(x,y)=
\begin{cases}
\sqrt{(\,|x|-0.8 \sin(\theta))^2+(y-0.8 \cos(\theta))^2}-0.2,
& \cos(\theta)|x|>\sin(\theta)y,\\[6pt]
\big|\sqrt{x^2+y^2}-0.8\big|-0.2,
& \cos(\theta)|x|\le \sin(\theta)y,
\end{cases}
\end{align*}
where $\theta = \pi/3$.
The solution remains close to the equilibrium state determined by the nonlinear reaction term throughout most of the interior of the domain $\Omega$, while the boundary condition enforces a rapid adjustment near $\partial\Omega$. 
This mismatch gives rise to boundary layers along the entire boundary, within which the solution undergoes a sharp transition from its interior level to the prescribed boundary value $u = 1$.

\begin{figure}[h]
\centering
\includegraphics[width=0.9\textwidth]{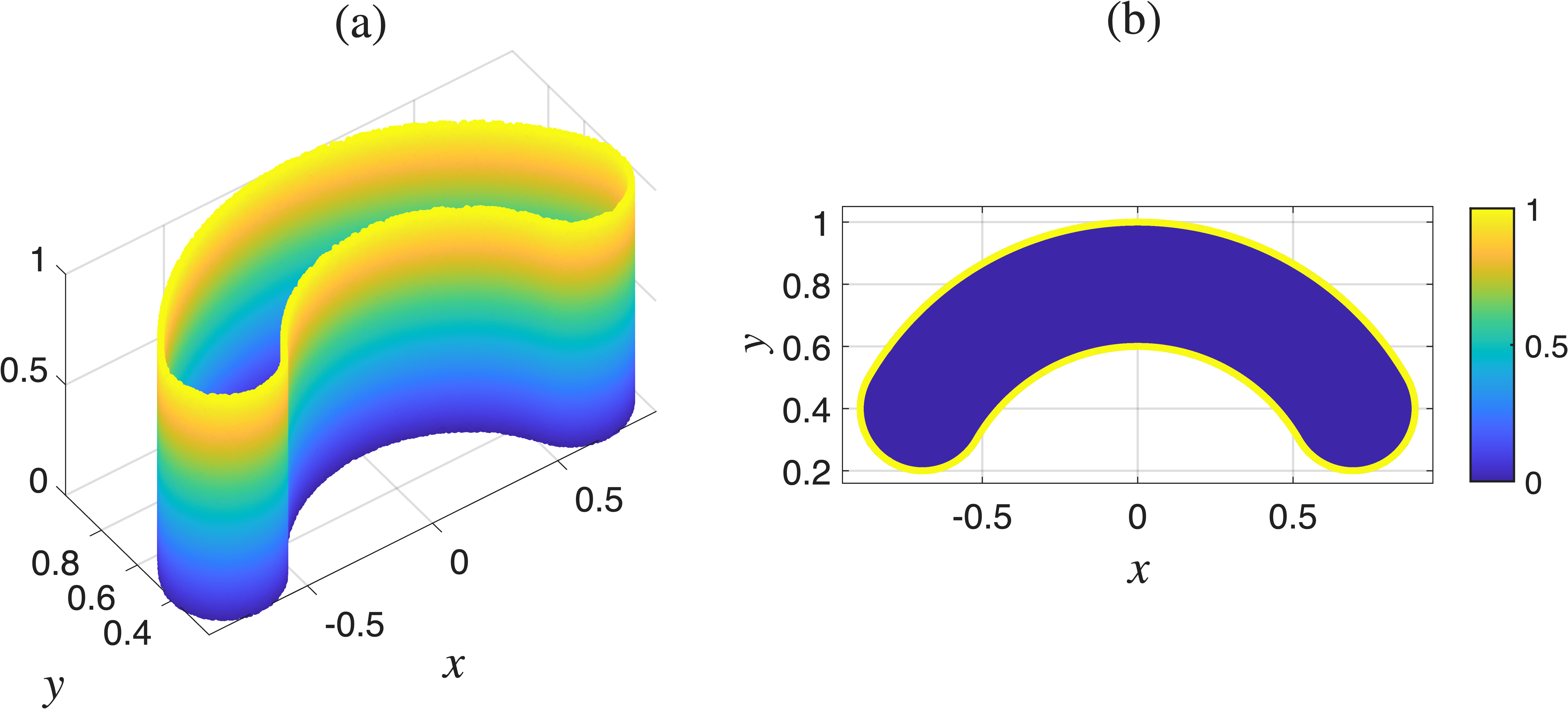}
\caption{Numerical solution of the two-dimensional Poisson--Boltzmann equation on an irregular domain with $\eps=10^{-10}$. 
  (a) Neural-network predicted solution, showing boundary layers formed along the entire boundary.
  (b) Top view of the solution, demonstrating that the solution remains close to zero in the interior region.}
\label{Fig:EX7}
\end{figure}

We employ the proposed neural network with 350 trainable parameters. 
Although an analytic solution is not available for this example, we conduct an experiment with $\eps = 10^{-10}$, for which the final loss value attains an order of $10^{-13}$. 
The numerically predicted solution, shown in Fig.~\ref{Fig:EX7} with both side and top views, exhibits the expected boundary-layer behavior. Moreover, the proposed framework can be directly applied to such nonlinear problems posed on irregular domains without any modification to the underlying network architecture.

\section{Conclusion and future work}\label{sec:conclusion}
In this work, we addressed the inherent challenges of learning solutions to singularly perturbed PDEs characterized by sharp boundary-layer structures. We proposed a unified learning framework featuring a weighted loss formulation that enables the automatic detection of  boundary-layer locations without relying on explicit asymptotic analysis. 
The proposed methodology directly adopts the standard PINN framework, consisting only of residual terms associated with the governing PDE and boundary conditions, thereby making the implementation straightforward.
Numerical experiments, on both scalar equations and systems of PDEs, demonstrated stable and robust performance for extremely small perturbation parameters $\eps = 10^{-10}$ while maintaining high solution accuracy across diverse domain geometries.

Beyond the specific class of singularly perturbed problems considered in this study, the present framework demonstrates that multiscale difficulties in neural-network-based PDE solvers need not be addressed through explicit asymptotic decompositions. Instead, carefully designed loss formulations can implicitly encode scale interactions and mitigate stiffness-induced optimization imbalance. This perspective suggests that the resolution of boundary-layer phenomena can be treated as an optimization-level challenge rather than solely a representation-level construction, offering a compelling alternative to classical perturbation-based strategies.

Several directions remain for future investigation. Extending the proposed framework to more general multiscale PDEs exhibiting both boundary and interior layers is a natural next step. Further analysis of the optimization dynamics induced by the weighted loss formulation could provide deeper insights into its robustness within stiff and strongly multiscale regimes. In addition, exploring the connections between the proposed optimization-driven mechanism and operator-learning architectures could offer a scalable approach for learning multiscale solution operators across varying parameter regimes. Finally, extending the framework to three-dimensional problems also presents significant challenges in terms of computational complexity and optimization stability, constitutes another important avenue for future research.

\section*{Acknowledgement}



W.-F. Hu, P.-W. Hsieh, and T.-S. Lin acknowledge support from the National Science and Technology Council (NSTC), Taiwan, under Grants NSTC 114-2628-M-008-002-MY4 (W.-F. Hu), NSTC 114-2115-M-005-001 (P.-W. Hsieh), and NSTC 114-2124-M-390-001 (W.-F. Hu and T.-S. Lin).


\bibliography{ref.bib}

%
%

\end{document}